\documentclass[]{article}
\usepackage[english]{babel}
\usepackage{multind}
\usepackage{lineno}
\usepackage{graphicx,multicol}
\usepackage{amssymb}
\usepackage{amsmath}
\newcommand{\mmk}[1]{\marginpar{\raggedright \small \em MK: #1}}

\newcommand{\qed}{\hfill$\diamond$\\\vspace{2mm}}
\newcommand{\pf}{{\bf Proof: }}
\newcommand{\be}{\begin{enumerate}}
\newcommand{\ee}{\end{enumerate}}
\newcommand{\bd}{\begin{description}}
\newcommand{\ed}{\end{description}}

\newcommand{\beq}{\begin{equation}}
\newcommand{\eeq}{\end{equation}}

\newtheorem{theorem}{Theorem}[section]
\newtheorem{lemma}[theorem]{Lemma}
\newtheorem{proposition}[theorem]{Proposition}

\newtheorem{corollary}[theorem]{Corollary}
\newtheorem{conjecture}[theorem]{Conjecture}
\newtheorem{question}[theorem]{Question}
\newtheorem{problem}[theorem]{Problem}

\newtheorem{claim}{Claim}

\def\p_1ed{\hbox{\kern1pt\vrule height6pt width4pt
depth1pt\kern1pt}\medskip}

\newcommand{\Z}{{\mathbb{Z}}}

\newcommand{\scri}{{\mathcal{I}}}

\newcommand{\scrp}{{\mathcal{P}}}
\newcommand{\scrq}{{\mathcal{Q}}}

\newcommand{\sm}{\setminus}

\newcommand{\UG}{U\!G}

\setpagewiselinenumbers
\def\pro{\noindent {\bf Proof. }}

\begin{document}
\bibliographystyle{plain}
\title{Antistrong digraphs
\thanks{This
    work was initiated while Bang-Jensen and Jackson attended the program ``Graphs,
    hypergraphs and computing'' at Institute Mittag-Leffler, spring
    2014. Support from IML is gratefully acknowledged}}
\author{J\o{}rgen Bang-Jensen\thanks{Department of Mathematics and
    Computer Science, University of Southern Denmark, Odense DK-5230,
    Denmark (email:jbj@imada.sdu.dk). The research of Bang-Jensen was
    also supported by the Danish research council under grant number
    1323-00178B} \and St\'ephane Bessy\thanks{LIRMM, Universit\'e
    Montpellier,161 rue Ada 34392 Montpellier Cedex 5 FRANCE
    (email:stephane.bessy@lirmm.fr)}\and Bill
  Jackson\thanks{School of Mathematical Sciences, Queen Mary University
    of London, Mile End Road London E1 4NS, England
    (email:b.jackson@qmul.ac.uk)} \and Matthias
  Kriesell\thanks{Institut f\"ur Mathematik, Technische Universit\"at
    Ilmenau, Weimarer Strasse 25, D--98693 Ilmenau, Germany
(email:matthias.kriesell@tu-ilmenau.de)}} \maketitle

\begin{abstract}
An antidirected trail in a digraph is a trail (a walk with no arc
repeated) in which the arcs alternate between forward and backward
arcs. An antidirected path is an antidirected trail where no vertex is
repeated. We show that it is NP-complete to decide whether two
vertices $x,y$ in a digraph are connected by an antidirected path,
while one can decide in linear time whether they are connected by an
antidirected trail. A digraph $D$ is antistrong if it contains an antidirected
$(x,y)$-trail starting and ending with a forward arc for every
choice of $x,y\in V(D)$. We show that antistrong connectivity can be
decided in linear time. We discuss relations between antistrong
connectivity and other properties of a digraph and show that the
arc-minimal antistrong spanning subgraphs of a digraph are the
bases of a matroid on its arc-set. We show that one can determine in
polynomial time the minimum number of new arcs whose addition to $D$
makes the resulting digraph the arc-disjoint union of $k$ antistrong
digraphs. In particular, we determine the minimum number of new arcs
which need to be added to a digraph to make it antistrong. We use
results from matroid theory to characterize graphs which have an
antistrong orientation and give a polynomial time algorithm for
constructing such an orientation when it exists. This immediately
gives analogous results for graphs which have a connected bipartite
2-detachment. Finally, we study arc-decompositions of
antistrong digraphs and pose several problems and conjectures.\\

\noindent{}{\bf Keywords:} antidirected path, bipartite
representation, matroid, detachment, anticonnected digraph
\end{abstract}

\section{Introduction}

We refer the reader to~\cite{bang2009} for notation and terminology
not explicitly defined in this paper. An {\bf antidirected path} in
a digraph $D$ is a path in which the arcs alternate between forward
and backward arcs. The digraph $D$ is said to be {\bf anticonnected}
if it contains an antidirected path between $x$ and $y$ for every
pair of distinct vertices $x,y$ of $D$. Anticonnected digraphs were
studied in~\cite{chartrandUM51}, where several properties such as
antihamiltonian connectivity have been considered.
We will show  in Theorem~\ref{antipathcheck} below that it is
NP-complete to decide whether a given digraph contains an antidirected
path between given vertices.

Our main purpose is to introduce a related connectivity property
based on the concept of a {\bf forward antidirected trail}, i.~e.~a
walk with no arc repeated which begins and ends with a forward arc
and in which the arcs alternate between forward and backward arcs. A
digraph $D$ is {\bf antistrong} if it has at least three vertices
and contains a forward antidirected $(x,y)$-trail for every pair of
distinct 
vertices $x,y$ of $D$. We say that
$D$ is {\bf $k$-arc-antistrong} if it has at least three vertices
and contains $k$ arc-disjoint forward antidirected $(x,y)$-trails
for all distinct $x,y \in V(D)$.

\medskip

The paper is organized as follows. First we show that, from an
algorithmic point of view, anticonnectivity is not an easy concept to
work with, since deciding whether a digraph contains an antidirected
path between a given pair of vertices is NP-complete. Then we move to
the main topic of the paper, antistrong connectivity, and show that
this relaxed version of anticonnectivity is easy to check
algorithmically. In fact, we show in Section~\ref{ASsec} that there is
a close relation between antistrong connectivity of a digraph $D$ and
its so called bipartite representation $B(D)$, namely $D$ is
antistrong if and only if $B(D)$ is connected. This allows us in
Section~\ref{ASaugsec} to find the minimum number of new arcs we need
to add to a digraph which is not antistrong so that the resulting
digraph is antistrong. Furthermore, using the bipartite representation
we show in Section~\ref{matroidsec} that the arc-minimal antistrong
spanning subdigraphs of a digraph $D$ form the bases of a matroid on
the arc-set of $D$.
More generally, we show that the subsets of $A$ which contain no
closed antidirected trails are the independent sets of a matroid on
$A$. In Section~\ref{ASorsec} we study the problem of deciding whether
a given undirected graph has an antistrong orientation. We show how to
reduce this problem to a matroid problem and give a characterization
of those graphs that have an antistrong orientation. For the convenience of readers who are not familiar with matroids, we also provide, as an appendix,  a purely graph theoretical proof  of the key step which is Lemma \ref{lem:mat3}. Both proofs can be converted
to  polynomial time algorithms which either finds an antistrong
orientation of the given input graph $G$ or produce a certificate
which shows that $G$ has no such orientation. In
Section~\ref{detachsec} we show that being orientable as an antistrong
digraph can be expressed in terms of connected 2-detachments of graphs
(every vertex $v$ is replaced by two copies $v',v''$ and every
original edge $uv$ becomes an edge between precisely one of the 4
possible pairs $u'v',u'v'',u''v',u''v''$) with the extra requirement
that the 2-detachment is bipartite and contains no edge of the form
$u'v'$ or $u''v''$. This imediately leads to a characterization of
graphs having such a 2-detachment.  Finally, in Section~\ref{nosepsec}
we show that one can decide in polynomial time whether a given digraph
$D$ has a spanning antistrong subdigraph $D'$ so that $D-A(D')$ is
connected in the underlying sense (while it is NP-hard to decide
whether a given digraph contains a non-separating strong spanning
subdigraph).\\ We conclude the paper with some remarks and open
problems.

\section{Anticonnectivity}\label{anticonsec}
It was shown in~\cite{chartrandUM51} that every connected graph $G$
has an anticonnected orientation. This can be seen by considering a
breath first search tree rooted at some vertex $r$. Let
$\{r\}=L_0,L_1,L_2,\ldots{},L_k$ be the distance classes of $G$.
Orient all edges between $r$ and $L_1$ from $r$ to these vertices,
orient all edges between $L_1\cup L_3$ and $L_2$ from $L_2$ to
$L_1\cup L_3$, orient all edges from $L_4$ to $L_3\cup L_5$ etc.
Finally, orient all the remaining, not yet oriented edges
arbitrarily.

We will show that it is NP-complete to decide if a digraph has an
antidirected path between two given vertices. We need the following
result which is not new, as it follows from a result in~\cite{gabowITSE2} on the vertex analogue,  but we include a new and short proof for
completeness.
\begin{theorem}
\label{avoidpairs}
It is NP-complete to decide for a given graph $G=(V,E)$, two specified
vertices $x,y\in V$ and pairs of distinct edges ${\cal
  P}=\{(e_1,f_1),(e_2,f_2),\ldots{},(e_p,f_p)\}$, all from $E$,
whether $G$ has an $(x,y)$-path which avoids at least one edge from
each pair in $\cal P$.
\end{theorem}

\pf We first slightly modify a very useful polynomial reduction, used in many papers such
as~\cite{bangTCS438}, from
3-SAT to a simple path problem and then show how to extend this to a reduction
from 3-SAT to the problem above. For simplicity our proof uses
multigraphs but it is easy to change to graphs.

Let $W[u,v,p,q]$ be the graph (the variable gadget) with vertices
$\{u,v,y_1,y_2,\dots{}y_p,z_1,z_2,\ldots z_q\}$ and the edges of the
two $(u,v)$-paths $uy_1y_2\ldots{}y_pv, uz_1z_2\ldots{}z_qv$.
\\
Let ${\cal F}$ be an instance of 3-SAT with variables
$x_1,x_2,\ldots{},x_n$ and clauses $C_1,C_2,\ldots{},C_m$.
The ordering of the clauses $C_1,C_2,\ldots{},C_m$ induces an ordering
of the occurrences of a variable $x$ and its negation $\overline{x}$ in
these. With each variable $x_i$ we associate a copy of
$W[u_i,v_i,p_i+1,q_i+1]$ where $x_i$ occurs $p_i$ times and
$\overline{x_i}$ occurs $q_i$ times in the clauses of $\cal F$. Identify
end vertices of these graphs by setting $v_i=u_{i+1}$ for
$i=1,2,\ldots{},n-1$. Let $s=u_1$ and $t=v_n$ and denote by $G'$ the
resulting graph. In $G'$ we respectively denote by $y_{i,j}$ and
$z_{i,j}$ the vertices $y_j$ and $z_j$ in the copy of $W$ associated with
the variable $x_i$.
\\
Next, for each clause $C_i$ we associate this with 3 edges from $G'$
as follows: assume $C_i$ contains variables $x_j,x_k,x_l$ (negated or
not). If $x_j$ is not negated in $C_i$ and this is the $r$th copy of
$x_j$ (in the order of the clauses that use $x_j$), then we associate
$C_i$ with the edge $y_{j,r}y_{j,r+1}$ and if $C_i$ contains
$\overline{x_j}$ and this is the $k$th occurrence of $\overline{x_j}$, then we
associate $C_i$ with the edge $z_{j,k}z_{j,k+1}$. We make similar
associations for the other two literals of $C_i$. Thus for each clause
$C_i$ we now have a set $E_i$ of three distinct edges
$e_{i,1},e_{i,2},e_{i,3}$ from $G'$ and $E_i\cap E_j=\emptyset$ for
$i\neq j$.
\\
Now it is easy to check that $G'$ has an $(s,t)$-path which avoids at
least one edge from each of the sets $E_1,E_2,\ldots{},E_m$ if and
only if $\cal F$ is satisfiable. Indeed, the $(s,t)$-path goes through the
`$z$-vertices' of the copy of $W$ associated with $x_i$ if
and only if  $x_i$ is
set to {\sc true} to satisfy $\cal F$.

Let us go back to the original problem. Let $H$ be the multigraph
consisting of vertices $c_0,c_1,\ldots{},c_m$ and three edges (denoted
$f_{i,1},f_{i,2},f_{i,3}$) from $c_{i-1}$ to $c_i$ for $i\in \{1,\dots,m\}$.
Let $G$ denote the multigraph we obtain from $G'$ and $H$ by identifying
$t$ and $c_0$. Let $x=s$ and $y=c_m$.  Finally, form three disjoint
pairs of arcs $(e_{i,1},f_{i,1}), (e_{i,2},f_{i,2}),(e_{i,3},f_{i,3})$
between $E_i$ and $\{f_{i,1},f_{i,2},f_{i,3}\}$ for every $i \in \{1 \dots m\}$.\\ By the observations
above it is easy to check that $G$ has an $(x,y)$-path which avoids
at least one arc from each of the forbidden pairs if and only if $\cal F$ is satisfiable. \qed

\begin{theorem}
\label{antipathcheck}
It is NP-complete to decide whether a given digraph contains an
antidirected path between given vertices $x,y$.
\end{theorem}

\pf The following proof is due to Anders Yeo (private communication,
April 2014).  Let $G=(V,E)$ be a graph with two specified vertices
$x,y\in V$ and pairs of distinct edges ${\cal
  P}=\{(e_1,f_1),(e_2,f_2),\ldots{},(e_p,f_p)\}$, all from $E$. We
will show how to construct a digraph $D_G$ with specified vertices
$s,t$ such that $D_G$ contains an antidirected $(s,t)$-path if and
only if $G$ has an $(x,y)$-path which avoids at least one edge from
each pair in $\cal P$. Since the construction can be done in
polynomial time this and Theorem~\ref{avoidpairs} will imply the
result.

Let $k$ be the maximum number of pairs in ${\cal P}$ involving the
same edge from $E$. Let $P$ be an antidirected path of length $2k+2$
which starts with a forward arc (and hence ends with a backward
arc). Now construct $D_G$ as follows: start from $G$ and first replace
every edge $uv$ with a private copy $P_{uv}$ of $P$ (no internal
vertices are common to two such paths). Then for each pair
$(e_i,f_i)\in {\cal P}$ we identify one sink of $P_{e_i}$ with one
source of $P_{f_i}$ so that the resulting vertex has in- and out-degree $2$.
By the choice of the length of $P$ we can identify in pairs, i.~e.~no three vertices will be identified.
Note that all the original vertices of $G$ will be sources in $D_G$.
The remaining (new vertices) will be called {\it internal} vertices.

Finally let $s=x$ and $t=y$. We claim that $D_G$ has an antidirected
$(s,t)$-path if and only if $G$ has an $(x,y)$-path which uses at most
one edge from each of the pairs in ${\cal P}$. Suppose first that
$xx_1x_2\ldots{}x_{r-1}x_ry$ is a path in $G$ which uses at most one
edge from each of the pairs in ${\cal P}$. Then
$P_{xx_1}P_{x_1x_2}\ldots{}P_{x_{r-1}x_r}P_{x_ry}$ is an antidirected
$(s,t)$-path in $D_G$ (no vertex is repeated since the identifications
above where only done for paths corresponding to pairs in ${\cal
  P}$). Conversely, suppose $D_G$ contains an antidirected
$(s,t)$-path $Q$.
By the way we identified vertex pairs according to ${\cal P}$, the internal
vertices have in- and out-degree at most $2$, and if an internal vertex
is on two paths $P_{e_i},P_{f_i}$ then it has both its in-neighbours
on $P_{e_i}$ and both its out-neighbours on $P_{f_i}$. This implies
that $Q$ will either completely traverse a path $P_{e_i}$ or not touch
any internal vertex of that path. Hence it cannot traverse both
$P_{e_i}$ and $P_{f_i}$ if $(e_i,f_i)\in {\cal P}$, and it follows that
if we delete all internal vertices of $Q$ and add back the edges of
$G$ corresponding to each of the traversed paths, we obtain an
$(x,y)$-path in $G$ that uses at most one edge from each pair in
${\cal P}$. \qed

\section{Properties of antistrong digraphs}\label{ASsec}
It follows from our definition that every pair of vertices of an antistrong digraph is joined by a trail of odd length. This immediately gives

\begin{lemma}
\label{nobipAS}
No bipartite digraph is antistrong.
\end{lemma}

For every digraph $D$ we can associate an undirected bipartite graph
which contains all the information we need to study antistrong
connectivity.  The {\bf bipartite representation}~\cite[Page
  19]{bang2009} of a digraph $D=(V,A)$ is the bipartite graph
$B(D)=(V'\cup V'',E)$, where $V'=\{v'|v\in V\}$, $V''=\{v''|v\in V\}$
and $E=\{v'w''|vw\in A\}$.

\begin{proposition}
\label{ASchar}
Let $D=(V,A)$ be a digraph with $|V|\ge 3$. The following are equivalent.
\begin{enumerate}
\item $D$ is antistrong
\item $B(D)$ is connected.
\item For every choice of distinct 
vertices $x,y$, the digraph $D$ contains both an
  antidirected $(x,y)$-trail $T_{x,y}$ of even length starting on a
  forward arc and an antidirected $(x,y)$-trail $\bar{T}_{x,y}$ of
  even length starting on a backward arc.
\end{enumerate}
\end{proposition}


\pf Suppose (a) holds. Then, following the edges corresponding to the
arcs of a forward antidirected $(x,y)$-trail, $B(D)$ contains an
$(x',y'')$-path for every pair of distinct vertices $x,y\in V$. Now,
if $x$ and $y$ are distinct vertices of $D$, we choose a third vertex
$z$ in $D$ ($z\neq x$ and $z\neq y$), and the union of an
$(x',z'')$-path and a $(z'',y')$-path contains an $(x',y')$-path in
$B(D)$. Similarly we obtain an $(x'',y'')$-path in $B(D)$ for every
pair of distinct vertices $x,y\in V$. Finally, for any $x\in V$, an
$(x',x'')$-path in $B(D)$ can be found in the union of an $(x',y')$-path
and a $(y',x'')$-path, where $y$ is a vertex of $D$ distinct from $x$.
Hence $(a)\Rightarrow{}(b)$ holds. Conversely, $(b)\Rightarrow{}(a)$ holds, since
any $(x',y'')$-path in $B(D)$ corresponds to a forward antidirected
$(x,y)$-path in $D$ which starts and ends with a forward arc.\\
Now to prove $(b)\Rightarrow{}(c)$, it suffices to remark that
$T_{x,y}$ and $\bar{T}_{x,y}$ correspond to an
$(x',y')$-path and an $(x'',y'')$-path in $B(D)$, respectively.
Finally, to see that $(c)\Rightarrow{}(b)$ holds, it suffices to show that if (c) holds, then
$B(D)$ contains an $(x',y'')$-path for all $x,y\in V$ (possibly
equal). This follows by considering a neighbour $z''$ of $x'$ and a
$(z'',y'')$-path in $B(D)$.  \qed

Proposition~\ref{ASchar} immediately implies the next result.
\begin{corollary}
\label{checkAS}
One can check in linear time whether a digraph is antistrong.
\end{corollary}





%

Recall that a digraph is {\bf $\mathbf{k}$-strong} if it has at least
$k+1$ vertices and it remains strong after deletion of any set of at
most $k-1$ vertices.  The digraph obtained from three
disjoint independent sets $X,Y,Z$ each of
size $k$ by adding all arcs from $X$ to $Y$, from $Y$ to $Z$, and from $Z$ to $X$ is $k$-strong. However,
$B(D)$ has three connected components. This shows that no condition
on the strong connectivity will
guarantee that a digraph is antistrong.\\

Recall that $D$ is $k$-arc-antistrong if it contains $k$ arc-disjoint
forward antidirected $(x,y)$-trails for every ordered pair of distinct
vertices $x,y$.  We can check in time $O(mk)$ whether a digraph has
$k$ arc-disjoint forward antidirected $(x,y)$-trails for given
vertices $x,y$, because they correspond to edge-disjoint $(x',y'')$-paths in
$B(D)$ whose existence can be checked by using flows, see
e.g.~\cite[Section 5.5]{bang2009}. So we can check in polynomial time
if a digraph is $k$-arc-antistrong.

\begin{theorem}
If $D$ is $2k$-arc-antistrong, then it contains
$k$ arc-disjoint antistrong spanning subdigraphs.
\end{theorem}
\pf Since $D$ is $2k$-arc-antistrong, $B(D)$ is $2k$-edge-connected.
We can now use Nash-Williams'~theorem
(see~\cite[Theorem~9.4.2]{bang2009} for instance) to deduce that
$B(D)$ has $k$ edge-disjoint spanning trees.
Proposition~\ref{ASchar} now gives the required  set of $k$
arc-disjoint antistrong spanning subdigraphs of $D$.
\qed

\begin{theorem}
There exists a polynomial time algorithm which for a given digraph $D$ and a
natural number $k$ either returns $k$ arc-disjoint spanning antistrong
subdigraphs
of $D$ or correctly answers that
no such set exists.
\end{theorem}

\pf This follows from the fact that such subdigraphs exist if and
only if $B(D)$ has $k$ edge-disjoint spanning trees, and the
existence of such trees can be checked via Edmonds' algorithm for
matroid partition~\cite{edmondsJRNBS69}. \qed

The corresponding problem for containing two arc-disjoint {\em
strong} spanning subdigraphs is NP-complete (see
e.g.~\cite[Theorem~13.10.1]{bang2009}).

\begin{theorem}
\label{twospstrongD} It is NP-complete to decide whether a digraph
$D$ contains two spanning strong subdigraphs $D_1,D_2$ which are
arc-disjoint.
\end{theorem}

\section{Antistrong connectivity augmentation}\label{ASaugsec}
Note that
every complete digraph on at least 3
vertices is antistrong. Hence it is natural to ask for the minimum
number of new arcs one has to add to a digraph in order to make it
antistrong.

\begin{theorem}
\label{augmenttoAS}
There exists a polynomial time algorithm for finding, for a given digraph
$D=(V,A)$ on at least 3 vertices, a minimum cardinality set of new arcs
$F$ such that the digraph $D'=(V,A \dot{\cup} F)$ is antistrong.
\end{theorem}
\pf Let $D$ be a digraph on $n\geq 3$ vertices which is not
antistrong. By Proposition~\ref{ASchar}, its bipartite representation
$B(D)$ is not connected. First observe that in the bipartite
representation each new arc added to $D$ will correspond to an arc
from a vertex $u'$ of $V'$ to a vertex $v''\in V''$ such that $u\neq
v$ back in $V$. So we are looking for the minimum number of new edges
of type $u'v''$ with $u\neq v$ whose addition to $B(D)$ makes it
connected while preserving the bipartition $V',V''$. Note that, as long as
$n\geq 3$, in which case $B(D)$ has at least 6 vertices, we can always
obtain a connected graph by adding edges that are legal according to
the definition above.  So the number of edges we need is exactly the
number of connected components of $B(D)$ minus one.\footnote{This
  number is also equal to $(2n-1)-r(A)$ where $r$ is the rank function
  of the matroid $M(D)$ which we define in Section~\ref{matroidsec}.}

To find an optimal augmentation we add all missing edges between
$V'$ and $V''$ to $B(D)$, except for those of the form $v'v''$ and
give the new edges cost 1, while all original edges get cost 0. Now find a
minimum weight spanning tree in the resulting weighted complete
bipartite graph. The edges of cost 1 correspond to an optimal
augmenting set back in $D$. \qed

The complexity of the analogous question for $k$-arc-antistrong
connectivity is open.\\

\begin{problem}
\label{addtogetkAS} Given a digraph $D$ and a natural number $k$,
can we find in polynomial time a minimum cardinality set of new arcs
whose addition to $D$ results in a digraph $D'$ which is
$k$-arc-antistrong?
\end{problem}

Problem~\ref{addtogetkAS} is easily seen to be equivalent to the
following problem on edge-connectivity augmentation of bipartite
graphs.

\begin{problem}
\label{bipaug}
Given a natural number $k$ and a bipartite graph $B=(X,Y,E)$ with
$|X|=|Y|=p$ which admits a perfect matching $M$ in its bipartite
complement, find a minimum cardinality set of new edges $F$ such that
$F\cap M=\emptyset$ and $B+F$ is $k$-edge-connected and bipartite with the
same bipartition as $B$.
\end{problem}


Theorem~\ref{augmenttoAS} can be extended to find the minimum number
of new arcs whose addition to $D$ gives a digraph with $k$
arc-disjoint antistrong spanning subdigraphs $D_1,\ldots{},D_k$,
provided that $V(D)$ is large enough to allow the existence of $k$
such subdigraphs. Note that since each $D_i$ needs at least $2n-1$
arcs and we do not allow parallel arcs, we need $n$ to be large
enough, in particular we must have $n\geq 2k+1$.

\begin{theorem}
\label{addtogetkdisjointAS}
There exists a polynomial time algorithm for determining, for a given
digraph $D$ on at least 3 vertices, whether one can add some edges to
$D$ such that the resulting digraph is simple (no parallel arcs) and
has $k$ arc-disjoint antistrong spanning subdigraphs. In the case when
such a set exists, the algorithm will return a minimum cardinality set
of arcs $A'$ such that $D'=(V,A\cup A')$ contains $k$ arc-disjoint
antistrong spanning subdigraphs.
\end{theorem}
\pf This follows from the fact that the minimum set of new arcs is
exactly the minimum number of new edges, not of the form $v'v''$ that
we have to add to $B(D)$ such that the resulting bipartite graph is simple and
has $k$ edge-disjoint spanning trees. This number can be found using
matroid techniques as follows. Add all missing edges from $V'$ to
$V''$ and give those of the form $v'v''$ very large cost (larger than
$2nk$) and the other new edges cost 1. Now, if the resulting complete
bipartite digraph $K_{n,n}$ has $k$-edge-disjoint spanning trees of
total cost less than $2nk$, then the set of new edges added will form
a minimum augmenting set and otherwise no solution exists.  Recall
from matroid theory that $k$ edge-disjoint spanning trees in $K_{n,n}$
correspond to $k$ edge-disjoint bases in the cycle matroid
$M(K_{n,n})$ of $K_{n,n}$ which again corresponds to an independent
set of size $k(2n-1)$ in the union $M=\bigvee_{i=1}^kM(K_{n,n})$.
This means that we
can solve the problem by finding a minimum cost base $B$ of $M$ and
then either return the arcs which correspond to edges of cost 1 in $B$
or decide that no solution exists when the cost of $B$ is more than
$2kn$. We leave the details to the reader.\qed

\section{A matroid for antistrong connectivity}\label{matroidsec}

Having seen the equivalence between antistrong connectivity of digraph
$D$ on $n$ vertices and connectivity of its bipartite representation
$B(D)$ (see Proposition~\ref{ASchar}),
and recalling from matroid theory that $B(D)$ is connected if
and only if the cycle matroid $M(B(D))$ has rank $|V(B(D))|-1$, it is
natural to ask how antistrong connectivity can be expressed as a
matroid property on $D$ itself.

For $F\subseteq A$, we denote by $h(F)$ and $t(F)$ the numbers of
vertices that are heads, respectively tails, of one or more arcs in
$F$.

Recall that the independent sets of the cycle matroid $M(G)$ of a
graph $G=(V,E)$ are
those subsets $I\subseteq E$ for which we have $|I'| \leq
\nu(I')-1$ for all $\emptyset\neq I'\subseteq I$, where $\nu(I')$
is the number of end vertices of the edges in $I'$. Inspired by this
we define set $I$ of arcs in a digraph $D=(V,A)$ to be {\bf
independent} if
\begin{equation}
\label{goodset} |I'| \leq h(I')+t(I')-1 \hspace{3mm}\mbox{for all
}\emptyset \neq I'\subseteq I,
\end{equation}
A set $S\subseteq A$ is {\bf dependent} if it is not independent.

\begin{proposition}
\label{goodisforest} Let $D=(V,A)$ be a digraph. A subset
$I\subseteq A$ is independent if and only if the corresponding edge
set $I$ in $B(D)$ forms a forest. Every inclusion-minimal dependent
set $S\subseteq A$ corresponds to a cycle in $B(D)$ and
conversely. 
\end{proposition}
\pf Suppose $I\subseteq A$ is independent and consider the
corresponding edge set $\tilde I$ in $B(D)$. If $\tilde I$ is not a
forest, then some subset $\tilde I'\subseteq \tilde I$ will be a
cycle $C$ in $B(D)$ with $p$ vertices in each of $V',V''$ for some
$p\geq 2$. The set $\tilde I'$ corresponds to a set $I'\subseteq I$
with $h(I')+t(I')-1=p+p-1<2p=|I'|$, contradicting that $I$ is
independent. The other direction follows from the fact that every
forest $F$ in $B(D)$ spans at least $|E(F)|+1$ vertices in $B(D)$
and every subset of a forest is again a forest. The last claim
follows from the fact that every minimal set of edges which does not
form a forest in $B(D)$ forms a cycle in $B(D)$. \qed


The previous proposition implies that a set of arcs of a digraph
is dependent if and only if it contains a closed trail of even
length consisting of alternating forward and backward arcs. We will
refer to such a trail as a {\bf closed antidirected trail}, or {\bf CAT}
for short.

\begin{theorem}\label{antimatroid}
Let $D=(V,A)$ be a digraph and  $\cal I$ be the family of all
independent sets of arcs in $D$. Then $M(D)=(A,{\cal I})$ is a
graphic matroid with rank equal to
the size of a largest collection of arcs containing no closed
alternating trail.
\end{theorem}
\pf It follows immediately from Proposition~\ref{goodisforest} that
a set $I$ belongs to ${\cal I}$ if and only if the corresponding
edge set $\tilde I$ is independent in the cycle matroid on $B(D)$.
\qed

\begin{theorem}
A digraph $D$ is antistrong if and only if $M(B(D))$ has rank $2|V|-1$.
\end{theorem}
\pf The rank of $M(B(D))$ equals the size of a largest acyclic set of
edges in $B(D)$. This has size $2|V|-1$ precisely when $B(D)$ has a
spanning tree $H$. Back in $D$, the arcs corresponding to $E(H)$
contain antidirected forward
trails between any pair of
distinct vertices. \qed

\section{Antistrong orientations of graphs}\label{ASorsec}

Recall that, by Robbins' theorem (see e.g.~\cite[Theorem
  1.6.1]{bang2009}) a graph $G$ has a strongly connected orientation
if and only if $G$ is 2-edge-connected. For antistrong orientations we
have the following consequence of Proposition~\ref{ASchar} which
implies that there is no lower bound to the (edge-) connectivity which guarantees
an antistrong orientation of a graph.

\begin{proposition}
\label{nobipOK}
No bipartite graph can be oriented as an antistrong digraph.
\end{proposition}

The purpose of this section is to
characterize graphs
which can be oriented as antistrong digraphs.

\begin{theorem}\label{orient}
Suppose $G=(V,E)$ and $|E|=2|V|-1$. Then $G$ has an antistrong
orientation if and only if
\begin{eqnarray}
|E(H)|&\leq& 2|V(H)|-1 \mbox{ for all nonempty subgraphs $H$ of $G$, and}\label{condition}\\
 |E(H)|&\leq& 2|V(H)|-2 \mbox{ for all nonempty bipartite subgraphs $H$ of $G$.}\label{bipcond}\\
\nonumber
\end{eqnarray}
\end{theorem}


We derive Theorem~\ref{orient} from the following characterization
of graphs which can be oriented as  digraphs with no closed
antidirected trail (CAT).

\begin{theorem}
\label{CATfree} A graph $G=(V,E)$ has an orientation with no CAT if
and only if $G$ satisfies (\ref{condition}) and (\ref{bipcond}).  In
particular no $n$ vertex graph with at least $2n$ edges and no $n$
vertex bipartite graph with at least $2n-1$ edges admits a CAT-free
orientation.
\end{theorem}

It is not hard to see that Theorem~\ref{CATfree} implies
Theorem~\ref{orient}. Assume that Theorem~\ref{CATfree} holds and
consider a graph $G=(V,E)$ with $|E|=2|V|-1$. Suppose that $G$ has an
antistrong orientation $D$. Then $B(D)$ is connected by
Proposition~\ref{ASchar}.  As $B(D)$ has $2|V|-1=|V(B(D))|-1$ edges it
is a tree. So $D$ is a CAT-free orientation of $G$ and, by
Theorem~\ref{CATfree}, conditions (\ref{condition}) and
(\ref{bipcond}) hold for $G$.  Conversely, if (\ref{condition}) and
(\ref{bipcond}) hold for $G$, then $G$ has a CAT-free orientation by
Theorem~\ref{CATfree}, and we can deduce as above that this
orientation is also an antistrong orientation of $G$.\\

We next show that (\ref{condition}) and (\ref{bipcond}) are necessary
conditions for a CAT-free orientation. For the necessity of
(\ref{condition}), suppose that some nonempty subgraph $H$ has $|E(H)|\geq
2|V(H)|$ and that $D$ is any orientation of $G$. Then $B(D)$ has at
least $2|V(H)|$ edges between $V(H)'$ and $V(H)''$, implying that it
contains a cycle. Hence $D$ is not CAT-free. The necessity of
(\ref{bipcond}) can be seen as follows. Suppose $H$ is a bipartite
subgraph on $2|V(H)|-1$ edges and let $\vec{H}$ be an arbitrary
orientation of $H$.
Since no bipartite graph has an antistrong orientation it
follows that $B(\vec{H})$ is not connected, and, as
it has $2|V(H)|-1=|V(B(\vec{H}))|-1$ edges, it
contains a cycle. This corresponds to a CAT in
$\vec{H}$.\\

Most of the remainder of this section is devoted to a proof of
sufficiency in Theorem~\ref{CATfree}. We first show that, for an
arbitrary graph $G'=(V',E')$, the edge sets of all subgraphs $G$ of
$G'$ which satisfy (\ref{condition}) and (\ref{bipcond}) are the
independent sets of a matroid on $E'$.  We then show that this matroid
is the matroid union of the cycle matroid and the `even bicircular
matroid' of $G'$ (defined below).  This allows us to partition the edge-set of a graph
$G$ which satisfies (\ref{condition}) and (\ref{bipcond}) into a
forest and an `odd pseudoforest'. We then use this partition to define
a CAT-free orientation of $G$.  We first recall some results from
matroid theory. We refer a reader unfamiliar with submodular functions
and matroids to~\cite{frank2011}.

Suppose $E$ is a set and $f:2^E\to \Z$ is a submodular, nondecreasing
set function which is nonnegative on $2^E \setminus
\{\emptyset\}$. Edmonds~\cite{E1970}, see~\cite[Theorem
  13.4.2]{frank2011}, showed that $f$ induces a matroid $M_f$ on $E$
in which $S\subseteq E$ is independent if $|S'|\leq f(S')$ for all
$\emptyset\neq S'\subseteq S$.
 The rank of a subset $S\subseteq E$ in $M_f$ is given by the min-max formula
 \begin{equation}\label{eq:minmax}
r_f(S)=\min_\scrp\left\{\left|S\sm\bigcup_{T\in \scrp}T\right|+\sum_{T\in\scrp}f(T)\right\},
\end{equation}
where the minimum is taken over all subpartitions $\scrp$ of $S$
(where a {\bf subpartition} of $S$ is a collection of pairwise
disjoint nonempty subsets of $S$). Note that the matroid $M(D)$ defined in
the previous section is induced on the arc-set of the digraph $D$ by
the set function $h+t-1$.

Given a graph $G=(V,E)$ and $S\subseteq E$ let $G[S]$ be the {\bf
subgraph induced} by $S$ i.e. the subgraph of $G$ with edge-set $S$
and vertex-set all vertices incident to $S$. Let $\nu,\beta:2^E\to
\Z$ by putting $\nu(S)$ equal to the number of vertices incident to
$S$, and $\beta(S)$ equal to the number of bipartite components of
$G[S]$. It is well known that $\nu$ is submodular, nondecreasing,
and nonnegative on $2^E$ and that $M_{\nu-1}(G)$ is the cycle
matroid of $G$. The function $\nu-\beta$ is also
submodular, nondecreasing, and nonnegative on $2^E $ since it is the rank function of the matroid on $E$ whose independent sets  are the edge sets of the {\bf
odd pseudoforests} of
  $G$, i.~e.~subgraphs in which each connected component contains at
most one cycle, and if such a cycle exists then it is odd, see \cite[Corollary 7D.3]{Z}.
We will refer to this matroid as the {\bf
even bicircular matroid} of $G$.

The above mentioned properties of $\nu$ and $\nu-\beta$ imply that
$2\nu-1-\beta$ is submodular, nondecreasing, and nonnegative on $2^E
\setminus \{\emptyset\}$. We will show that the independent sets in
$M_{2\nu-1-\beta}(G)$ are the edge sets of the subgraphs which
satisfy (\ref{condition}) and (\ref{bipcond}).

\begin{lemma}\label{lem:mat1}
Let $G=(V,E)$ be a graph and $\scri=\{I\subseteq E\,:\mbox{ $G[I]$
  satisfies (\ref{condition}) and (\ref{bipcond})}\}$. Then $\scri$ is
the family of independent sets of the matroid
$M_{2\nu-1-\beta}(G)$. In addition, the rank of a subset $S\subseteq
E$ in this matroid is
$r_{2\nu-1-\beta}(S)=\min_\scrp\left\{|S\sm\bigcup_{T\in
  \scrp}T|+\sum_{T\in \scrp}\left(2\nu(T)-1-\beta(T)\right)\right\}$
where the minimum is taken over all subpartitions $\scrp$ of $S$.
\end{lemma}
\pro We first suppose that some $S\subseteq E$ is not independent in
$M_{2\nu-1-\beta}(G)$. Then we may choose a nonempty $S'\subseteq S$
with $|S'|>2\nu(S')-1-\beta(S')$, and subject to this condition, such
that $|S'|$ is as small as possible. The minimality of $S'$ implies
that $H=G[S']$ is connected. So $\beta(S')=1$ if and only if $H$ is
bipartite (and 0 otherwise) and we may now deduce that that
$H\subseteq G[S]$ fails to satisfy (\ref{condition}) or
(\ref{bipcond}).

We next suppose that $G[S]$ fails to satisfy (\ref{condition}) or
(\ref{bipcond}) for some $S\subseteq E$. Then there exists a nonempty
subgraph $H$ of $G[S]$ such that either $|E(H)|>2|V(H)|-1$, or $H$
is bipartite and $|E(H)|>2|V(H)|-2$. Then $S'=E(H)$ satisfies
$|S'|>2\nu(S')-1-\beta(S')$ so $S$ is not independent in
$M_{2\nu-1-\beta}(G)$.

The expression for the rank function of $M_{2\nu-1-\beta}(G)$
follows immediately from (\ref{eq:minmax}). \qed

The {\bf matroid union} of two matroids $M_1=(E,\scri_1)$ and
$M_2=(E,\scri_2)$ on the same ground set $E$ is the matroid $M_1\vee
M_2=(E,\scri)$ where $\scri=\{I_1\cup I_2\,:\,I_1\in \scri_1\mbox{
and } I_1\in \scri_1\}$. Suppose $f_1,f_2:E\to \Z$ are submodular,
nondecreasing, and nonnegative on $2^E \setminus \{\emptyset\}$.
Then $f_1+f_2$ will also be submodular, nondecreasing, and
nonnegative on $2^E \setminus \{\emptyset\}$ and hence will induce
the matroid $M_{f_1+f_2}$. Every independent set in $M_{f_1}\vee
M_{f_2}$ is independent in $M_{f_1+f_2}$, but the converse does not
hold in general. Katoh and Tanigawa~\cite[Lemma 2.2]{KT} have shown
that the equality $M_{f_1+f_2}=M_{f_1}\vee M_{f_2}$ does hold
whenever the minimum in formula (\ref{eq:minmax}) for the ranks
$r_{f_1}(S)$ and $r_{f_2}(S)$ is attained for the same subpartition
of $S$, for all $S\subset E$. This allows us to deduce

\begin{lemma} \label{lem:mat2}
For any graph $G=(V,E)$, we have $M_{2\nu-1-\beta}(G)=M_{\nu-1}(G)\vee
M_{\nu-\beta}(G)$.
\end{lemma}
\pro This follows from the above mentioned result of Katoh and
Tanigawa, and the facts that
$r_{\nu-1}(S)=\sum_{T\in\scrp}(\nu(T)-1)$  and
$r_{\nu-\beta}(S)=\sum_{T\in\scrp}(\nu(T)-\beta(T))$ where $\scrp$
is the partition of $S$ given by the connected components of $G[S]$
(since $r_{\nu-1}(S)$ and $r_{\nu-\beta}(S)$ are equal to the number
of edges in  a maximum forest and a maximum odd pseudoforest,
respectively, in $G[S]$). \qed

Lemma~\ref{lem:mat1} and Lemma~\ref{lem:mat2} immediately give the following.

\begin{lemma}\label{lem:mat3}
Let $G=(V,E)$ be a graph. Then $G$ satisfies
(\ref{condition}) and (\ref{bipcond}) if and only if $E$ can be
partitioned into a forest and an odd pseudoforest.
\end{lemma}

\noindent We provide an alternative graph theoretic proof of this
lemma in the appendix.

We next show that every graph whose edge set can can be partitioned
into a spanning tree and an odd pseudoforest has a CAT-free
orientation.\\


\begin{figure}[h]
\begin{center}
\includegraphics[height=85mm]{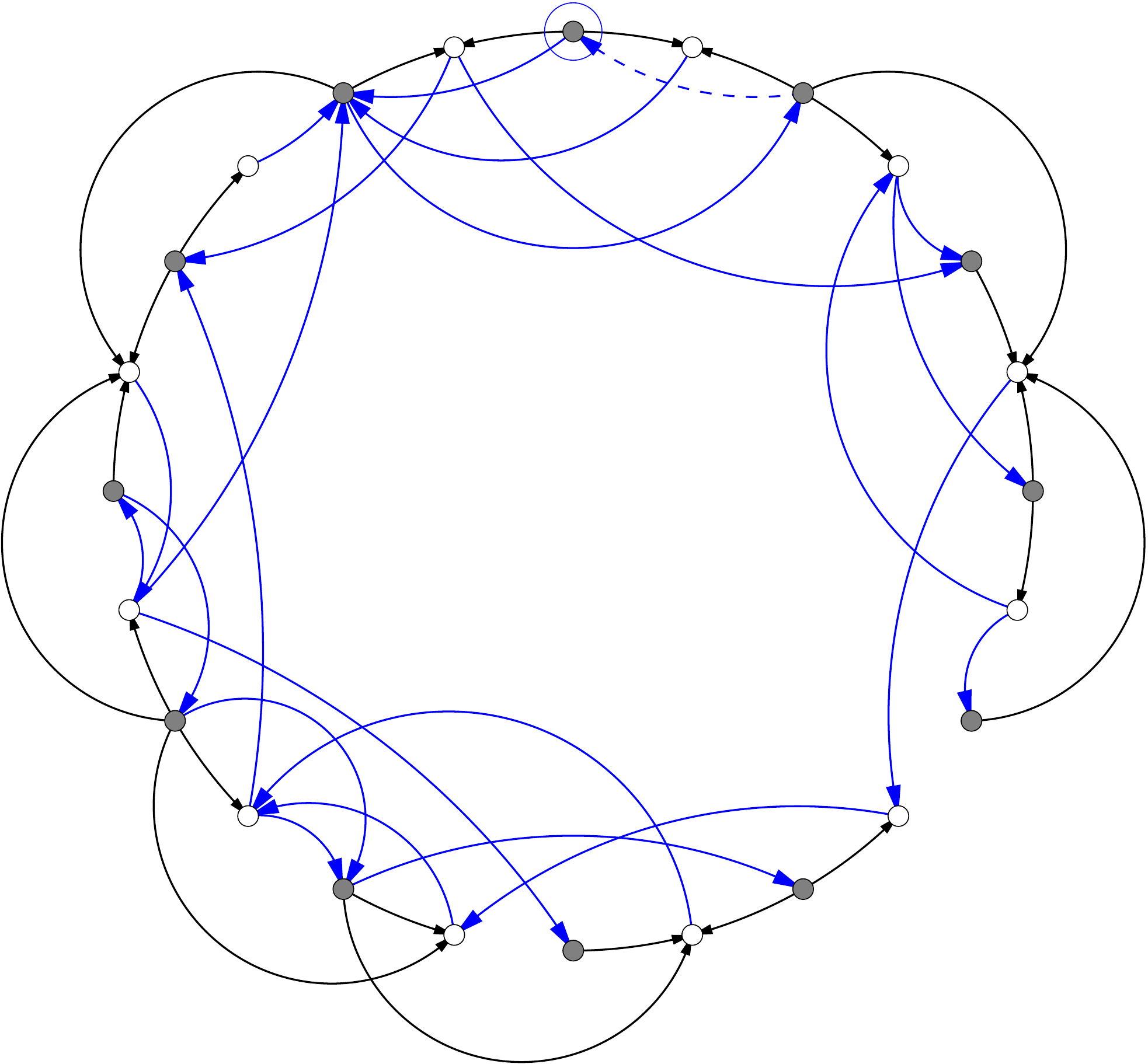}%
\end{center}
\caption{A CAT-free orientation of the union of a spanning tree $T$
  and a spanning odd pseudoforest $P$; $T$ governs the bipartition $X,Y$
  (white/grey), its edges are drawn outside (or on) the disk spanned by the
  vertices.  The edges of $P$ are embedded in the interior of that
  disk, the root vertex is the encircled topmost one, the precious
  edge is the dashed one.}
\label{twotreesfig}
\end{figure}
\begin{theorem}
\label{tree+forest}
Let $G$ be the edge-disjoint union of a spanning tree $T$ and an odd
pseudoforest $P$.\footnote{Note that $G$ may have parallel edges, but
  no more than two copies of any edge, in which case one copy is in
  $T$ and the other in $P$} Then $G$ has a CAT-free orientation. In
addition, such an orientation can be constructed in linear time given
$T$ and $P$.
\end{theorem}
\pf Let $X,Y$ be the unique (up to renaming the two sets)
bipartition of $T$ and orient all edges of $T$ from $Y$ to $X$. If
$P$ has no edges we are done since there are no cycles in $G$. Let
$P_1,\ldots{},P_k$ be the connected components of $P$. We shall show
that we can orient the edges of $P_1,\ldots{},P_k$ in such a way
that none of the resulting arcs of these (now oriented)
pseudoforests
$\vec{P_1},\ldots{},\vec{P_k}$
can belong to a closed antidirected trail. Clearly this will imply
the lemma. For each $P_i$ we choose a root vertex $r_i$ of $P_i$ as
follows. If $P_i$ is a tree then we choose $r_i$ to be an arbitrary
vertex of $P_i$. If $P_i$ contains an odd cycle $C_i$ then we choose
$r_i$ to be a vertex of $C_i$ such that $r_i$ has as many neighbours on $C_i$ as possible
in the same set of the
bipartition $(X,Y)$ as $r_i$. Since $C_i$ is odd we may choose one such neighbour
$s_i$ of $r_i$.
We will refer to the edge $r_is_i$ as a {\bf precious edge} of
$P_i$. Put $T_i=P_i-r_is_i$ if $P_i$ contains a cycle and otherwise
put $T_i=P_i$.

We orient the edges of $T_i$
as follows. Every edge of $P_i$ with one end in $X$ and the
other in $Y$ is oriented from $X$ to $Y$. Every edge $uv$ of $T_i$ with $u,v\in
X$ is oriented towards $r_i$ in $T_i$
(so if $v$ is closer to $r_i$ than $u$ in $T_i$ we orient the edge
from $u$ to $v$ and otherwise we orient it from $v$ to $u$, see
Figure~\ref{twotreesfig}). Every edge $pq$ of $T_i$ with $p,q\in Y$ is
oriented away from $r_i$ in $T_i$.
Finally, if $P_i$ contains a precious edge $r_is_i$, then we orient
$r_is_i$ from $r_i$ to $s_i$ if $r_i,s_i\in X$, and from $s_i$ to
$r_i$ if $r_i,s_i\in Y$.  Let $D$ denote the resulting orientation of
$G$.  The digraph $D$ can be constructed in linear time if we traverse
each tree $T_i$ by a breath first search rooted at $r_i$.

We use induction on $|E(P)|$ to show that the above construction results in a CAT-free digraph. As noted
above, this is true for the base case when $E(P)=\emptyset$. Suppose
that $E(P)\neq \emptyset$ and choose an edge $uv$ in some $P_i$
according to the following criteria. If $P_i$ is not a cycle then
choose $v$ to be a vertex of degree one in $P_i$ distinct from $r_i$
and $u$ to be the neighbour of $v$ in $P_i$. If $P_i$ is an odd cycle
then choose $v=r_i$ and $u=s_i$. We will show that $uv$ belongs to no
CAT in $D$.  By symmetry, we may suppose that $v\in X$.

We first consider the case when $v$ is a vertex of degree one in
$P_i$. Below $d^+(v),d^-(v)$ denote the out-degree, respectively, the
in-degree of the vertex $v$ in $D$.
We have two possible subcases:
\begin{itemize}
\item $u\in Y$. Since $v\in X$, we oriented $uv$ from $v$ to $u$. All
  the other edges incident to $v$ belong to $T$ and were oriented
  towards $v$. Then $d^+(v)=1$ and the arc $vu$ cannot be part of a
  CAT.
\item $u\in X$. Since $v\in X$, we oriented $uv$ from $v$ to $u$ (as $u$
  is closer to $r_i$ than $v$ in $T_i$). As previously we have $d^+(v)=1$ and
  the arc $vu$ cannot be part of a CAT.
\end{itemize}
Since $D-uv$ is CAT-free by induction, $D$ is also CAT-free.

We next consider the case when $P_i$ is an odd circuit.
Let $t_i$ be the neighbour of $r_i$ in $P_i$ distinct from $s_i$.
We again have  two possible subcases:
\begin{itemize}
\item $t_i\in X$.  Since $r_i\in X$, we oriented the edge $t_ir_i$ from $t_i$ to
  $r_i$. Then $d^+(r_i)=1$, and the arc $r_is_i$ cannot be part of a CAT. Since $D-r_is_i$ is CAT-free by induction, $D$ is also CAT-free.
\item $t_i\in Y$. Let $q_i$ be the neighbour of $s_i$ in $P_i$ which
  is distinct from $r_i$. The choice of $r_i$ implies that $q_i\in Y$,
  and hence that $s_iq_i$ is oriented from $s_i$ to $q_i$. Then
  $d^+(s_i)=1$, and the arc $s_iq_i$ cannot be part of a CAT. Since $D-s_iq_i$ can be obtained by applying our construction to $T\cup (P-s_iq_i)$ by choosing $s_i$ as the root vertex in $P_i-s_iq_i$ rather than $r_i$, it is CAT-free by induction. Hence $D$ is also CAT-free.
\end{itemize}

\qed

\noindent{}{\bf Proof of Theorem~\ref{CATfree} (sufficiency):} Let
$G=(V,E)$ be a graph satisfying (\ref{condition}) and
(\ref{bipcond}). By Lemma~\ref{lem:mat3}, $E$ can be partitioned into
a forest $F$ and an odd pseudoforest $P$. By adding a suitable set of
edges to $G$, we may assume that $|E|=2|V|-1$. (This follows by
considering the matroid $M_{2\nu-1-\beta}(2K_n)$ on the edge set of
the graph $2K_n$ with vertex set $V$ in which all pairs of vertices
are joined by two parallel edges. It is easy to check that $2K_n$ has
an edge-disjoint forest and odd pseudoforest with a total of $2|V|-1$
edges. Thus the rank of $M_{2\nu-1-\beta}(2K_n)$ is $2|V|-1$. Since
$E$ is an independent set in $M_{2\nu-1-\beta}(2K_n)$, it can be
extended to an independent set with $2|V|-1$ edges.) The fact that
$|E|=2|V|-1$ implies that $F$ is a spanning tree of $G$. We can now
apply Theorem~\ref{tree+forest} to deduce that $G$ has a CAT-free
orientation.  \qed

We have seen that Theorem~\ref{CATfree} implies Theorem~\ref{orient},
and hence that a graph $G=(V,E)$ has an antistrong orientation if and
only if the rank of $M_{2\nu-1-\beta}(G)$ is equal to $2|V|-1$. We can
now apply the rank formula (\ref{eq:minmax}) to characterize graphs
which admit an antistrong orientation.

\begin{theorem}\label{thm:char}
A graph $G=(V,E)$ has an antistrong orientation if and only if
\begin{equation}\label{ASorchar}
e(\scrq)\geq |\scrq|-1+b(\scrq)
\end{equation}
for all partitions $\scrq$ of $V$, where $e(\scrq)$ denotes the
number of edges of $G$ between the different parts of $\scrq$ and
$b(\scrq)$ the number of parts of $\scrq$ which induce bipartite
subgraphs of $G$.
\end{theorem}

\noindent
\pf Suppose that $G$ has no antistrong orientation. Then the rank of
$M_{2\nu-1-\beta}(G)$ is less than $2|V|-1$ so there exists a
subpartition $\scrp$ of $E$ such that
\begin{equation}\label{ASorchar1}
\alpha(\scrp):=\left|E\sm\bigcup_{T\in
\scrp}T\right|+\sum_{T\in\scrp}\left(2\nu(T)-1-\beta(T)\right)<
2|V|-1
\end{equation}
by Lemma \ref{lem:mat1}. We may assume that $\scrp$ has been chosen such that:
\begin{enumerate}
\item[(i)] $\alpha(\scrp)$ is as small as possible;
\item[(ii)] subject to (i), $|\scrp|$ is as small as possible;
\item[(iii)] subject to (i) and (ii), $|\bigcup_{T\in
  \scrp}T|$ is as large as possible.
  \end{enumerate}
   Let
$\scrp=\{E_1,E_2,\ldots,E_t\}$ and let $H_i=(V_i,E_i)$ be the
subgraph of $G$ induced by $E_i$ for all $1\leq i\leq t$.
We will show that $H_i$ is a (vertex-)induced connected subgraph of
$G$ and that $V_i\cap V_j=\emptyset$ for all $i\neq j$.

First, suppose that $H_i$ is disconnected for some $1\leq i\leq t$.
Then we have $H_i=H_i'\cup H_i''$ for some subgraphs
$H_i'=(V_i',E_i')$ and $H_i''=(V_i'',E_i'')$ with $V_i'\cap
V_i''=\emptyset$. Let $\scrp'=(\scrp\sm\{E_i\})\cup\{E_i',E_i''\}$. We
have
$$2\nu(E_i)-1-\beta(E_i)>2\nu(E_i')-1-\beta(E_i')+2\nu(E_i'')-1-\beta(E_i'')$$ since
$\nu(E_i)=\nu(E_i')+\nu(E_i'')$ and $\beta(E_i)=
\beta(E_i')+\beta(E_i'')$. This implies that
$\alpha(\scrp')<\alpha(\scrp)$ and contradicts (i). Hence $H_i$ is
connected and $\beta(E_i)\in \{0,1\}$ for all $1\leq i \leq t$.

Next, suppose that $V_i\cap V_j\neq \emptyset$ for some $1\leq
i<j\leq t$. Let $\scrp'=(\scrp\sm\{E_i,E_j\})\cup\{E_i\cup E_j\}$. We
have
$$2\nu(E_i)-1-\beta(E_i)+2\nu(E_j)-1-\beta(E_j)\geq 2\nu(E_i\cup
E_j)-1-\beta(E_i\cup E_j)$$ since, if $|V_i\cap V_j|=1$, then
$\nu(E_i)+\nu(E_j)= \nu(E_i\cup E_j)+1$ and
$\beta(E_i)+\beta(E_j)\leq \beta(E_i\cup E_j)+1$,
and, if $|V_i\cap V_j|\ge 2$, then $\nu(E_i)+\nu(E_j)\geq
\nu(E_i\cup E_j)+2$ and $\beta(E_i)+\beta(E_j)\leq 2$.
This implies that $\alpha(\scrp')\leq\alpha(\scrp)$. Since
$|\scrp'|<|\scrp|$ this contradicts (i) or (ii). Hence $V_i\cap
V_j=\emptyset$ for all $i\neq j$.

Finally, suppose that $H_i\neq G[V_i]$. Then some $e\in
E\sm\bigcup_{T\in
  \scrp}T$ has both end vertices in $E_i$. Let $E_i'=E_i+e$ and
$\scrp'=\scrp - E_i +E_i'$. This implies that
$\alpha(\scrp')\leq\alpha(\scrp)$. Since $|\scrp'|=|\scrp|$ and
$|\bigcup_{T\in
  \scrp'}T|>|\bigcup_{T\in
  \scrp}T|$,
this contradicts (i) or (iii). Hence $H_i=G[V_i]$.

Let $\scrq$ be the partition of $V$ obtained from
$\{V_1,V_2,\ldots,V_t\}$ by adding the remaining vertices of $G$ as
singletons. Then $\left|E\sm\bigcup_{T\in \scrp}T\right|=e(\scrq)$
and
$\sum_{T\in\scrp}\left(2\nu(T)-1-\beta(T)\right)=2|V|-|\scrq|-b(\scrq)$.
We can now use (\ref{ASorchar1}) to deduce that $e(\scrq)<
|\scrq|-1+b(\scrq)$.

\smallskip

For the converse, suppose that $e(\scrq)< |\scrq|-1+b(\scrq)$ for
some partition $\scrq=\{V_1,V_2,\ldots,V_s\}$ of $V$. We may assume that $G[V_i]$ is connected for all $1\leq i \leq s$ since we can replace $V_i$ by the vertex sets of the components of $G[V_i]$ in $\scrq$ and maintain this inequality. Let
$G[V_i]=(V_i,E_i)$ for $1\leq i\leq s$ and  $\scrp=\{E_i\,:\,E_i\neq
\emptyset,\,1\leq i\leq s\}$.
 Then $\left|E\sm\bigcup_{T\in
\scrp}T\right|=e(\scrq)$ and
$\sum_{T\in\scrp}\left(2\nu(T)-1-\beta(T)\right)=2|V|-|\scrq|-b(\scrq)$. (Note that a set $V_i\in \scrq$ with $|V_i|=1$ has no corresponding edge set in $\scrp$ and contributes $2-1-1$ to the right hand side of the last equation.)
This implies that
$$\left|E\sm\bigcup_{T\in
\scrp}T\right|+\sum_{T\in\scrp}\left(2\nu(T)-1-\beta(T)\right)=
e(\scrq)+\sum_{T\in\scrp}\left(2\nu(T)-1-\beta(T)\right)<
2|V|-1$$ and hence $G$ has no antistrong orientation. \qed
\\[2mm]

\begin{corollary} Every $4$-edge-connected nonbipartite graph has an
antistrong orientation.
\end{corollary}
\pf Suppose $G=(V,E)$ is $4$-edge-connected and not bipartite and
let $\scrq$ be a partition of $V$. If $\scrq=\{V\}$ then
$e(\scrq)=0=|\scrq|-1+b(\scrq)$ since $G$ is not bipartite, and if
$\scrq\neq \{V\}$ then $e(\scrq)\geq 2|\scrq|\geq
|\scrq|-1+b(\scrq)$ since $G$ is $4$-edge-connected. Hence $G$ has
an antistrong orientation by Theorem~\ref{thm:char}. \qed

\begin{corollary}
\label{cor:3trees-antistrong}
Every nonbipartite graph with three edge-disjoint spanning trees has
an antistrong orientation.
\end{corollary}
\pf We give two proofs of this corollary.

 Suppose $G=(V,E)$ is a
nonbipartite graph with three edge-disjoint spanning trees and let
$\scrq$ be a partition of $V$.  If $\scrq=\{V\}$ then
$e(\scrq)=0=|\scrq|-1+b(\scrq)$ since $G$ is not bipartite, and if
$\scrq\neq \{V\}$ then $e(\scrq)\geq 3(|\scrq|-1)$ since $G$ has
three edge-disjoint spanning trees.
Since $|\scrq|\ge 2$,
$2(|\scrq|-1)\ge |\scrq| \ge b(\scrq)$ and $e(\scrq)\geq
|\scrq|-1+b(\scrq)$. Hence $G$ has an antistrong orientation by
Theorem~\ref{thm:char}.

 We could also remark that if $T_1$, $T_2$
and $T_3$ denote three edge-disjoint spanning trees of $G$, then there
exists $e\in G$ such that $T_1+e$ is not bipartite. Then depending if
$e\in T_2$ or not, $\{T_1+e,T_3\}$ or $\{T_1+e,T_2\}$ is an
edge-disjoint pair of a spanning odd pseudo-tree and a spanning tree
of $G$. Let $H$ denote this subgraph of $G$. Then using
Theorem~\ref{tree+forest}, $H$ has a CAT-free orientation which is also
an antistrong orientation of $H$ since $|E(H)|=2|V(H)|-1$. So $G$ has also an
antistrong orientation.
\qed

Corollary \ref{cor:3trees-antistrong} is tight in the sense that
there exist graphs with many edge-disjoint trees, two spanning and
the others missing just three vertices, which have no antistrong
orientation. Consider the graph $G$ obtained by identifying one
vertex of a complete bipartite graph $K_{k,k}$ and a complete graph
$K_4$. Then $G$ has no antistrong orientation. Indeed, consider the
partition $\scrq$ of $V(G)$ into four parts: the copy of $K_{k,k}$,
and one part for each remaining vertex of $K_4$.  We have
$e(\scrq)=6 < |\scrq|-1+b(\scrq)= 4-1+4$ and then $G$ has no
antistrong orientation by Theorem~\ref{thm:char}.


\medskip

Since $M_{2\nu-1-\beta}(G)=M_{\nu-1}(G)\vee M_{\nu-\beta}(G)$, we
can use Edmonds' matroid partition algorithm~\cite{edmondsJRNBS69}
to determine the rank of $M_{2\nu-1-\beta}(G)$ in polynomial time,
and hence determine whether $G$ has an antistrong orientation.
Moreover, when such an orientation exists, we can use Edmonds'
algorithm to construct an edge-disjoint spanning tree and
pseudoforest with a total of $2|V|-1$ edges, and then use the
construction from the proof of Theorem~\ref{tree+forest} to obtain
the desired antistrong orientation in polynomial time. This gives
\begin{corollary}
There exists a polynomial time algorithm which finds, for a given
input graph $G$, either an antistrong orientation $D$ of $G$, or a
certificate, in terms of a subpartition $\cal P$ which violates
(\ref{ASorchar}), that $G$ has no such orientation.
\end{corollary}

\section{Connected bipartite 2-detachments of graphs}\label{detachsec}

We now show a connection between antistrong orientations of a graph
$G$ and so-called detachments of $G$. We need only the special case of
2-detachments (see e.g.~\cite{nashwilliamsJLMS31} for results on
detachments).

A {\bf 2-detachment} of a graph $G=(V,E)$ is any graph $H=(V'\cup
V'',E')$
which can be obtained from $G$
by replacing every vertex $v\in V$ with two new vertices $v',v''$
and then for each original edge $uv$ adding precisely one of the
four edges $u'v',u'v'',u''v',u''v''$ to $E'$.

\begin{lemma}\label{detach}
A graph $G=(V,E)$ has an antistrong orientation if and only if $G$ has
a 2-detachment $H=(V'\cup V'',E')$ which is connected and bipartite
with bipartition $V',V''$ (we call such a 2-detachment {\bf good}).
\end{lemma}
\pf Suppose $G$ has a good 2-detachment $H=(V'\cup V'',E')$. Then
there are no edges of the form $u'v'$ and no edges of the form
$u''v''$. Hence the orientation $D$ that we get by orienting the edges
of the form $u'v''$ from $u$ to $v$ will be an antistrong orientation
of $G$ by Proposition~\ref{ASchar}.  Conversely, if $D$ is an
antistrong orientation of $G$, then $B(D)$ is a good 2-detachment of
$G$. \qed

We can now use Theorem~\ref{thm:char} and the subsequent remark to deduce the following.
\begin{theorem}
\label{good2detach} A graph $G=(V,E)$  has a good 2-detachment if
and only if
\begin{equation}
\label{detachcond}
e(\scrq)\geq |\scrq|-1+b(\scrq)
\end{equation}
for all partitions $\scrq$ of $V$. Furthermore, there exists a
polynomial time algorithm which returns such a $2$-detachment when it
exists and otherwise returns a certificate, in terms of a partition violating (\ref{detachcond}), that no such detachment exists.
\end{theorem}

\section{Non-separating antistrong spanning subdigraphs}\label{nosepsec}

While there are polynomial time algorithms for checking the existence of
two edge-disjoint spanning trees~\cite{edmondsJRNBS69}, or two
arc-disjoint branchings (spanning out-trees) in a digraph (see
e.g.~\cite[Corollary 9.3.2]{bang2009}), checking whether we can delete
a strong spanning subdigraph and still have a connected digraph is
difficult. Let $\UG(D)$ denote the underlying undirected graph of a digraph $D$.

\begin{theorem}\cite{bangTCS438}
\label{nondisconstrong}
It is NP-complete to decide whether a given digraph $D$ contains a
spanning strong subdigraph $H$ such that $\UG(D-A(H))$ is connected.
\end{theorem}

If we replace ``strong'' by ``antistrong'' above, the problem becomes
solvable in polynomial time.

\begin{theorem}
\label{deleteanticon}
We can decide in polynomial time for a given digraph $D=(V,A)$ on $n$
vertices whether $D$ contains a spanning antistrong subdigraph
$H=(V,A')$ such that $\UG(D-A')$ is connected.
\end{theorem}
\pf We may assume that $D$ is antistrong, since this can be checked in
linear time by verifying that $B(D)$ is connected.  Let $M_1=(A,{\cal
  I})$ be the cycle matroid of of the underlying graph $\UG(D)$ of $D$ and let $M_2=M(D)=(A,{\cal
  I}(D))$ be the matroid from Section~\ref{matroidsec} whose bases are
the antistrong sets consisting of $2n-1$ arcs. Let $M=M_1\vee{}M_2$ be
the union of the matroids $M_1,M_2$, that is, a set $X$ of arcs is
independent in $M$ if and only we can partition $X$ into $X_1, X_2$
such that $X_i$ is independent in $M_i$.
For each of the
matroids $M_1, M_2$ we can check in polynomial time whether a given
subset $X$ of arcs is independent in $M_1$ and $M_2$ (for $M_1$ we
need to check that there is no cycle in $\UG(D)[X]$ and for $M_2$ we
need to check that there is no cycle in the subgraph of $B(D)[E_X]$
induced by the edges $E_X$ corresponding to $X$ in $B(D)$%
). Thus it follows from
Edmonds' algorithm for matroid partitioning~\cite{edmondsJRNBS69} that
we can find a base of $M$ in polynomial time using the independence
oracles of $M_1,M_2$. The desired digraph $H$ exists if and only if
the size of a base in $M$ is $(2n-1)+(n-1)=3n-2$. \qed

A similar proof gives the following.

\begin{theorem}
We can decide in polynomial time whether a digraph $D$ contains
$k+\ell$ arc-disjoint spanning subdigraphs
$D_1,\ldots{},D_{k+\ell}$ such that
$D_1,\ldots{},D_k$ are antistrong and $\UG(D_{k+1}),\ldots{},\UG(D_{k+\ell})$
are connected.
\end{theorem}

\section{Remarks and open problems}

We saw in Theorem~\ref{twospstrongD} that it is NP-complete to decide whether a
given digraph contains two arc-disjoint spanning strong subdigraphs.  We
would be interested to know what happens if we modify the problem as
follows.

\begin{question}
Can we decide in polynomial time whether a digraph $D$ contains arc-disjoint
spanning subdigraphs $D_1,D_2$ such that $D_1$ is antistrong and $D_2$
is strongly connected?
\end{question}

Inspired by Theorem~\ref{deleteanticon} it is natural to ask the
following intermediate question.

\begin{question}
\label{anti+2con}
Can we decide in polynomial time whether a digraph $D$ contains arc-disjoint
spanning subdigraphs $D_1,D_2$ such that $D_1$ is antistrong and
$\UG(D_2)$ is 2-edge-connected?
\end{question}

The following conjecture was raised in~\cite{bangC24}.

\begin{conjecture}\cite{bangC24}
There exists a natural number $k$ such that every $k$-arc-strong
digraph has arc-disjoint strong spanning subdigraphs $D_1,D_2$.
\end{conjecture}

\noindent
Perhaps the following special case may be easier to study.

\begin{conjecture}
There exists a natural number $k$ such that every digraph $D$ which is
both $k$-arc-strong and $k$-arc-antistrong has arc-disjoint strong
spanning subdigraphs $D_1,D_2$.
\end{conjecture}

\begin{problem}
Does there exist a polynomial algorithm for deciding whether a given undirected graph $G$ has an orientation $D$ which is both strong and antistrong?
\end{problem}

\bigskip

\centerline{*}

\bigskip

{\bf Acknowledgement.}
Bang-Jensen and Jackson wish to thank Jan van den Heuvel for
stimulating discussions about antistrong connectivity.  They also
thank the Mittag-Leffler Institute for providing an excellent working
environment.


\pagebreak
\appendix
\section{Appendix: a graph theoretical proof of Lemma~6.6}

In this appendix, we give a graph theoretical proof of Lemma~6.6, recalled
below.\\

\noindent
{\bf Lemma~6.6} {\it
Let $G=(V,E)$ be a graph. Then $G$ satisfies}
\setcounter{equation}{1}
\begin{align}
|E(H)|\ &\leq \  2|V(H)|-1 \mbox{ for all nonempty subgraphs $H$ of $G$, and}\\
|E(H)|\ &\leq \  2|V(H)|-2 \mbox{ for all nonempty bipartite subgraphs $H$ of $G$}
\end{align}
{\it if and only if $E$ can be partitioned into a forest and an odd
pseudoforest.}\\

\pf
Recall that a {\bf pseudoforest} is a graph in which each connected
component contains at most one cycle, and it is called {\bf odd} if it does not contain even cycles.
A theorem due to Whiteley~\cite{whiteleyCM197} (see also~\cite{frank2011} p.367 for a
short proof based on Edmonds' branching theorem) asserts that a graph
satisfies condition~(2) if and only if its edge set can be
partitioned into a forest and a pseudoforest.
So let us denote by a {\bf 2-decomposition} of a graph $G=(V,E)$ a pair $(G_b,G_r)$
of spanning subgraphs $G_b=(V,E_b)$ and $G_r=(V,E_r)$ such that
$\{E_b,E_r\}$ is a partition of $E$ and $G_b$ is a forest of $G$ and $G_r$ is a pseudoforest of $G$.  We will call any sub-structure
--- edge, component, subgraph etc. --- of $G_r$ or of $G_b$ {\bf red} or {\bf black}, respectively,
and for a subgraph $H$ of $G$ we denote by $H_r$ and $H_b$ the subgraph of
$H$ induced by its red or black edges, respectively.

Without loss of generality we may assume that $G$ is connected, and
that $G_b$ is a spanning tree of $G$ (otherwise we could move
edges from $G_r$ to $G_b$ to make $G_b$ connected).
The {\bf canonical   bipartition} of a $2$-decomposition $(G_b,G_r)$ of a connected graph
$G$ is the unique bipartition given by any $2$-colouring of $G_b$.
Moreover, an edge of $G_r$ which does not cross this bipartition is called (as previously)
a {\bf precious edge} in $(G_b,G_r)$.  A 2-decomposition of a
2-decomposable graph is {\bf nice} if every red cycle contains at
least one precious edge.\\ First we establish the next claim.

\begin{claim}
\label{nicedecomp}
A connected graph which has a 2-decomposition admits a nice
2-decomposition if and only if
(\ref{bipcond}) holds.
\end{claim}
\noindent
\pf
First observe that for any subgraph $H$ of $G$ with at least one black edge, we have
$|E(H)|=|E(H_b)|+|E(H_r)|\le (|V(H_b)|-1) + |V(H_r)|= 2|V(H)|-1$.
For any red subgraph $H$ with at least one edge, we get $|E(H)|=|E(H_r)| \leq |V(H_r)| \leq 2|V(H)|-2$.
In particular (\ref{condition}) holds for every 2-decomposable graph.

The necessity is quite clear. Indeed, consider a nice
2-decomposition $(G_b,G_r)$ of $G$ and assume that (\ref{bipcond}) does not
hold. Thus there exists $H$ a bipartite subgraph of $G$ with
$|E(H)|=2|V(H)|-1$. So equality holds in the previous computation and we
have $|E(H_b)|=|V(H_b)|-1$ and $|E(H_r)|=|V(H_r)|$. In particular $H_b$ is a
spanning tree of $H$ and $H_r$ contains at least one cycle $C$.  As
$(G_b,G_r)$ is nice, $C$ contains a precious edge $xy$.  As $H_b$ is
connected, there exists a black path $P$ from $x$ to $y$ and $P$ has
even length because $x$ and $y$ belong to the same part of the
canonical bipartition of $(G_b,G_r)$. So $P\cup xy$ forms an odd cycle
of $H$, a contradiction.

Now let us prove the sufficiency. Let $(G_b,G_r)$ be a 2-decomposition
of $G$. A red component $R$ of the decomposition is {\bf bad} if it is not a tree
and its (hence unique) cycle does not contain any precious edges. If we remove from a bad
component $R$ all its precious edges, we obtain several connected
components, one of which contains the cycle of $R$. We call this
component the {\bf core} of $R$ and denote it by $c(R)$.  For
convenience we use $c(R)$ below to denote both a vertex set and the
red subgraph induced by these vertices. Note that $G[c(R)]$ is
bipartite as $c(R)$ contains no precious edge.
A {\bf sequence} of the decomposition
$(G_b,G_r)$ is a list ${\cal R}=(c(R_1),\dots ,c(R_i))$ of the cores of its bad red components in
decreasing order of cardinality.

Among all the 2-decompositions of $G$, we choose one whose sequences
${\cal R}= (c(R_1), \dots , c(R_i))$ satisfy

$(a)$ $i$ is minimal, and

$(b)$ subject to $(a)$, $|c(R_i)|$ is minimal.

\noindent
We will prove that this 2-decomposition $(G_b,G_r)$ is nice, that is,
${\cal R}=\emptyset$. Assume that it is not the case and consider
$\{X_1, \dots , X_p\}$ the black components of $G[c(R_i)]$ (that is,
the connected components of $G_b[c(R_i)]$). If $p=1$, then $G[c(R_i)]$
is connected in black, and as $G_r[c(R_i)]$ is unicyclic, the
bipartite graph $G[c(R_i)]$ violates $(\ref{bipcond})$, a
contradiction.  So we must have $p\ge 2$.
Now denote by $W_1,\dots,W_q$ the connected components of $G_b\setminus c(R_i)$ and construct
a graph $T'$ on $\{X_1, \dots ,X_p,W_1,\dots ,W_q\}$ by connecting two
vertices of $T'$ if there exists an edge in $G_b$ between the
corresponding components. In other words, we contract the (connected)
vertex sets $X_1, \dots ,X_p,W_1,\dots ,W_q$ in $G_b$ to single vertices.
So $T'$ is a tree. Finally we consider $T$ the minimal subtree of $T'$ containing
$\{X_1, \dots ,X_p\}$.  By definition the leaves of $T$ are in $\{X_1,
\dots ,X_p\}$ and as $p\ge 2$, $T$ has at least two such leaves.  So
we consider a leaf $X_k$ of $T$ which does not contain entirely the
red cycle of $R_i$ (this could occur even without violating
$(\ref{bipcond})$ if $X_k$ is not connected in red for instance). We
denote by $W_{k'}$ the only neighbour of $X_k$ in $T$.
Now, we specify two edges, one black and one red in order to `change their color' and
obtain a contradiction.  First denote by $uv$ the unique black edge
between $X_k$ and $W_{k'}$. We suppose that $u\in X_k$ and $v\in
W_{k'}$ (so $v\notin c(R_i)$). Now we look at a 1-orientation of
$c(R_i)$ (this is an orientation of $c(R_i)$ in which every vertex has out-degree
at most 1) and consider a maximal oriented red path leaving $u$ with all
its vertices in $X_k$. As $X_k$ does not contain entirely the red
cycle of $R_i$, this path ends at a vertex $x\in X_k$ which has a
unique red out-neighbour $y\in X_{k''}$ for some $k''\neq k$. We
select this red edge $xy$.\\

Notice that the unique black path from
$x$ to $y$ contains the edge $uv$.  Indeed the unique path from $X_k$
to $X_{k''}$ in $T$ corresponds to the unique black path $P$ from
$X_k$ to $X_{k''}$ in $G$. As $X_k$ is a leaf of $T$, the first edge
of $P$ is $uv$ and as $X_k$ and $X_{k''}$ are respectively connected
in black, the unique black path from $x$ to $y$ contains $P$ and so it
contains the edge $uv$. This implies that
$(G_b+xy)- uv$ is also a spanning tree of
$G$. Moreover, its bipartition is the same as the bipartition of
$G_b$. Indeed as $xy$ is an edge lying inside the core of $R_i$, it is
not precious and $G_b+xy$ is still bipartite and has the same
bipartition as $G_b$. Removing $uv$ does not affect the bipartition
(because $(G_b+xy)- uv$ is connected).

To conclude, we focus on the red part of the new 2-decomposition
$((G_b+xy)- uv, (G_r+uv)-xy)$. By construction, the component $X$ of
$G_r-xy$ containing $u$ (and also $x$) is a red tree. Remark that
$X=R_i$ if and only if  $xy$ is an edge of the cycle of $R_i$.  If
$v$ does not belong to $R_i$, then by adding $uv$ we attached $X$ to
another red component in $(G_r+uv)-xy)$. As $X$ contains at least
one vertex, namely $u$, $|c(R_i)|$ has decreased, a contradiction to
$(b)$ in the choice of $(G_b,G_r)$ (or to $(a)$ if $X=R_i$). If $v$
belongs to $R_i$ but $v$ does not belong to $X$ (in this case we
have $X\neq R_i$), then $v\in R_i\setminus c(R_i)$ and by adding
$uv$ we attached $X$ to a vertex of $R_i\setminus c(R_i)$. Once
again, $|c(R_i)|$ has decreased, a contradiction to $(b)$ in the
choice of $(G_b,G_r)$. Finally if $v$ belongs to $X$ then adding
$uv$ produces a new red unicyclic component.  However as the red path
in $G_r$ from $v$ to $u$ starts in $R_i\setminus c(R_i)$ and ends in
$c(R_i)$, it contains a precious edge.  So that newly created red unicyclic
component is not bad, and $|c(R_i)|$ has decreased.
Hence, again, we either contradict $(b)$, or $(a)$ if $X=R_i$.
~\\

Now to finish the proof, we will show how to go from a nice
2-decomposition of a connected graph to a decomposition into a
spanning tree and an odd pseudoforest (i.e. a pseudoforest in which
every cycle has odd length). Let $G=(V,E)$ be a connected graph which
admits a nice 2-decomposition and consider a nice 2-decomposition
$(G_b,G_r)$ of $G$ with a minimum number of even red cycles. We will
show by contradiction that this decomposition has no even red cycle.
Assume it is not the case and denote by $C_1, \dots ,C_l$ the even red
cycles of $G_r$. In each of these, select a precious edge
$e_i=x_iy_i$ and let  $X=\{x_1,y_1,x_2,y_2,\ldots{},x_l,y_l\}$.
Exchanging two edges between $G_b$ and
$G_r$ will modify the structure of $G_b$, and some previously selected
precious edges could become not precious any more. To avoid this we will
find a vertex $u$ with the following property\\

${\cal P}:$ There exists a
component $B$ of $G_b\setminus u$ such that one of the following hold:
\begin{itemize}
\item
 $B\cap X$ contains
only one element and this  is not in the same component of $G_r$ as $u$ ({\it Case A}).
\item $(B\cup \{u\})\cap X$ contains exactly two elements
and they are the end vertices of some $e_i$ ({\it Case B}).
\end{itemize}

First assume we have found such a vertex $u$ and let us see how to
conclude, depending of which the two cases A or B we are in.\\

\noindent{}{\it Case A}. Denote by $x_i$ the only element of $B\cap
X$ and by $B_i$ the red component of $G_r$ containing $x_i$.  As $u$
does not belong to $B_i$, we can find an edge $vw$ along the black
path in $G_b$ from $u$ to $x_i$ such that $w$ and $x_i$ are in the
same component of $G_b-vw$, $w$ belongs to $B_i$ and $v$ does not
belong to $B_i$. So we exchange the colors of $vw$ and $x_iy_i$.  The
graph $B_i-x_iy_i$ is a tree and when we add $vw$ to $G_r$ we connect
this tree to another component of $G_r$.  The component of $G_b-vw$
containing $v$ is a tree containing all the vertices of $X$ except
$x_i$.  So the precious edges $e_j$ with $j\neq i$ are still precious
edges, and this is also the case in $(G_b-vw)+x_iy_i$ which is a
spanning tree of $G$. So, we reduce the number of even red cycle of
the nice 2-decomposition $(G_b,G_r)$, a contradiction.\\

\noindent{}{\it Case
  B}. Denote by $x_i$ and $y_i$ the two elements of $(B\cup \{u\})\cap
X$ and also by $B_i$ the red component of $G_r$ containing the
precious edge $x_iy_i$. If the black path $P$ in $G_b$ between $x_i$
and $y_i$ is not totally contained in $B_i$ then we can select a
vertex $u'$ not belonging to $B_i$ along this path and end up in the
previous case with $u$ replaced by $u'$. So $P$ is totally contained in $B_i$.
 Then, as
$P+x_iy_i$ is an odd cycle (because $x_iy_i$ is precious), we can find
along $P$ two consecutive vertices $vw$ which are in the same part of
the bipartition induced by the bipartite graph $G_r[B_i]$. So we
exchange $x_iy_i$ and $vw$. As previously $G_b+x_iy_i-vw$ is a
spanning tree of $G$ such that all the edges $e_j$ with $j\neq i$ are
still precious and $vw$ is also precious.  The graph
$G_r-x_iy_i+vw$ is now a pseudoforest,
and we have reduced the number of even red cycles of the nice
2-decomposition $(G_b,G_r)$, a contradiction.

Finally, let us see how to find  a vertex $u$ in $G$ which has property $\cal P$. Consider
$T'$ the minimal subtree of $G_b$ containing all the vertices of the
set $X$. In particular all the leaves of $T'$ are elements of
$X$. Then build the tree $T$ from $T'$ by replacing iteratively each
vertex of degree 2 in $T'$ and not belonging to $X$ by an edge linking
its two neighbours in $T'$. The vertices of $T$ are now vertices of
$X$ or have degree at least three in $T$. Assume first that a leaf $f$
of $T$ has its neighbour $f'$ in $X$. Denote by $B$ the component of
$G_b\setminus f'$ containing $f$. By construction $f$ is the unique
element of $B\cap X$. We select $u=f'$. If $f$ and $f'$ are in
different components of $G_r$ then we are in {\it Case A},
otherwise we are in {\it Case B}.\\
So we can
assume that all the leaves of $T$ are neighbours of vertices of $T$
which are not in $X$ and have degree at least three in $T$. Consider
now a leaf $f'$ of the tree obtained from $T$ by removing its
leaves. Denote by $L$ the set of neighbours of $f'$ in $T$ which are
leaves of $T$.
If $|L|=2$ and $L$ consists of the end vertices of some $e_i$ then we choose $u=f'$ and are in {\it Case B}.
Otherwise, let $B_i$ be the component of $G_r$ containing $f'$ and
consider a vertex $f$ of $L$ not belonging to $B_i$ (this exists as $|B_i \cap X|=2$). Then we choose
$u=f'$ and $B$ to be the component of $G_b\setminus f'$ containing $f$
and we are in {\it Case A}. \qed

\end{document}